\documentclass[1p,times]{elsarticle} 
\usepackage{amsmath,amssymb,amsthm}         
\usepackage{mathtext}
\usepackage{epstopdf}
\usepackage{multicol}
\usepackage[noend]{algorithmic}
\usepackage{pgfplots}
\usepackage[colorinlistoftodos,textwidth=4cm,shadow]{todonotes}
\usepackage[colorlinks,urlcolor=blue]{hyperref}
\usepackage{algorithmic,algorithm}

\newtheorem{remark}{Remark}[section]

\newcounter{Ivan}

\newcounter{Maxim}

\begin{document}
\begin{frontmatter}
\title{Grid-based electronic structure calculations: the tensor decomposition approach}
\author[skoltech]{M.V.~Rakhuba} 
\ead{rakhuba.m@gmail.com}

\author[skoltech,inm]{I.V.~Oseledets}
\ead{i.oseledets@skoltech.ru}

\address[skoltech]{Skolkovo Institute of Science and Technology, Novaya St. 100, 143025 Skolkovo, Moscow Region, Russia.}

\address[inm]{Institute of Numerical Mathematics, Russian Academy of Sciences. Gubkina St. 8, 119333 Moscow, Russia.}

\begin{abstract}
We present a fully grid-based approach for solving Hartree-Fock and all-electron Kohn-Sham equations based on low-rank approximation of three-dimensional electron orbitals.
Due to the low-rank structure the total complexity of the algorithm depends linearly with respect to the one-dimensional grid size. 
Linear complexity allows for the usage of fine grids, e.g. $8192^3$ and, thus, cheap extrapolation procedure.

We test the proposed approach on closed-shell atoms up to the argon, several molecules and clusters of hydrogen atoms.
All tests show systematical convergence with the required accuracy.
\end{abstract}
\end{frontmatter}

\section{Introduction}

Electronic structure calculations arise in a variety of applications such as condensed matter physics or  material and drug design.
In this paper we focus on the development of efficient and accurate numerical techniques for solving Hartree-Fock (HF) \cite{hartree-1957} and Kohn-Sham (KS) equations \cite{ks-1965}. 
From computational point of view HF and KS equations are nonlinear eigenvalue problems with three-dimensional integro-differential operators. 
A standard way to solve these equations is to approximate the solution in a subspace of globally supported basis functions, e.g. gaussian functions.
This is a classic topic and a lot of software packages are available.
The choice of these basis functions is determined by the complexity of the the iterative process and the desired accuracy. 
It also introduces the \emph{basis set error}, which might be difficult to control.
Methods based on a sequence of embedded subspaces allow for rigorous error control. Among these approaches 
we should mention multiresolution approach based on wavelet function representation \cite{har-multiresdft-2004}, finite element methods on unstructured grids \cite{gav-femdft-2013}, finite difference method \cite{saad-fdm-1994} and  the projector-augmented wave method \cite{mjbw-gpaw-2005}.
Standard finite element or finite difference method on uniform grids for a 3D-dimensional problem has $\mathcal{O}(n^3 \log n)$ complexity (using Fast Fourier Transform for the evaluation of the integral operators). On non-uniform meshes we can not use FFT directly, thus the complexity can reach $\mathcal{O}(n^6)$, where $n$ is the one-dimensional grid size.

One of the promising ways to reduce the complexity, proposed and studied in details in the papers by Khoromskij and Khoromskaia~\cite{vkhs-2el-2013,vekh-hartree-2008,mpi-chem3d-2009,venera-phd,khor-chem-2011,vekh-qtt-hartree-2011,vekh-blackbox-2013, khor-ml-2009, kk-qsurvey-2015}  is to use \emph{tensor decomposition} of the vectors of coefficients. 
It was shown that in many cases, the tensor of coefficients can be well-approximated 
by the \emph{Tucker decomposition}, and the storage complexity goes down to $\mathcal{O}(n \log n)$, where $n$ is one-dimensional grid-size. This allows for huge uniform grids to be used. However, to compute the solution, they still employ an intermediate basis set (which can be arbitrary), and successively compute the Fock matrix and update the coefficients with respect to this basis. The main difficulty in solving the KS/HF equation directly in the tensor format is that standard iterative methods should be modified in order to keep the solution in the tensor format, and it requires the development of a new iterative method. 

Our main contribution in this paper is an efficient $\mathcal{O}(n \log n)$ black-box solver which is based on low-rank tensor decompositions, and it does not require any additional basis sets. Thus, it can be used to control the basis set error in HF/KS computations. We have developed an iterative process that solves the original problem while all intermediate three-dimensional arrays of coefficients are stored in the Tucker format.
Our method has several important components: we use the block integral iteration, fast low-rank convolution and derivative-free update of the Fock matrix. 
Finally, the extrapolation over several grid sizes is done, giving up to $\mathcal{O}(h^4)$ accuracy, where $h$ is the grid step.

We organize our paper as follows. In section \ref{sec:equations} we formulate the closed-shell Hartree-Fock and Kohn-Sham equations. 
In section \ref{sec:tensors} we discuss tensor decomposition approach and  the Tucker format.
Section \ref{sec:iterative} contains formulation of the Block-Green iteration and the derivative-free Fock matrix update formula.
Sections \ref{sec:discr} and \ref{sec:operations} are concerned with obtaining the discrete formulation and all operations involved in the tensor formats.
In section \ref{sec:complexity} we discuss the complexity of the proposed algorithm.
Finally, in section \ref{sec:num} we present numerical experiments on  closed-shell atoms up to argon, several molecules and cluster of hydrogen atoms, which illustrate sublinear rank dependence.

\section{Hartree-Fock and Kohn-Sham equations}\label{sec:equations}
Consider closed-shell electron system with $N_e$ electrons and $N = {N_e}/{2}$ orbitals, $N_{nuc}$ nuclei with charges $Z_{\alpha}$ located at $\bold{R}_{\alpha}, \alpha = 1, \ldots, N_{nuc}$.
Then the HF/KS equations can be written as \cite{jens-introqchem-2013}
\begin{equation}\label{hf:general_form}
H\left( \Phi \right) \phi_i \equiv \left( -\frac{1}{2} \Delta + V(\Phi)  \right) \phi_i = \lambda_i \phi_i, \quad i = \overline{1, N}
\end{equation}
where $\phi_i : \mathbb{R}^3\to \mathbb{R}^3$ denotes unknown orbitals that additionally satisfy the orthogonality constraints
$$
\int_{\mathbb{R}^3} \phi_i (\bold{r}) \phi_j (\bold{r}) \, d \bold{r}= \delta_{ij}.
$$
and $\lambda_i$ denotes orbital energies.
In the Kohn-Sham model
\begin{equation}\label{pot_general}
V(\Phi) \equiv \tilde V(\rho) = V_{ext} + V_{coul}(\rho) + V_{xc}(\rho),
\end{equation}
where 
$$
\rho(\bold{r}) = 2 \sum_{i=1}^{N} |\phi_i(\bold{r})|^2,
$$
is the electron density. 
The external potential $V_{ext}$ describes interaction between electrons and nuclei:
\begin{equation} \label{ext-pot}
V_{ext}(\bold{r}) = -\sum_{\alpha=1}^{N_{nuc}} \frac{Z_\alpha}{|\bold{r}-\bold{R}_\alpha |}.
\end{equation}
The potential $V_{coul}(\bold{r})$ is given as 
$$
V_{coul}(\bold{r}) = \int_{\mathbb{R}^3} \frac{\rho(\bold{r'})}{| \bold{r} - \bold{r'} |}\, d\bold{r'},
$$
and $V_{xc}$ depends only on the density and is responsible for exchange and correlation parts of the operator.
For the sake of simplicity we consider local density approximation (LDA) with Perdew and Zunger type functional \cite{pz-lda-1981}. We note that the concept can be extended to local spin density approximation (LSDA) or more accurate functionals.

In the HF equation the potential $V(\Phi)$ has the form
$$
V(\Phi) = V_{ext} + V_{coul}(\rho)  - \hat K(\Phi),
$$
 where
$$
\hat K(\Phi)\ \phi_i =\int_{\mathbb{R}^3} \frac{\tau(\bold{r}, \bold{r'})}{| \bold{r} - \bold{r'} |} \,\phi_i(\bold{r'}) \, d\bold{r'}
\quad\text{with}\quad 
\tau(\bold{r}, \bold{r'}) =  \sum_{i=1}^{N} \phi_i(\bold{r})\phi_i(\bold{r'}). 
$$
As $\hat K$ depends not only on density $\rho(\bold{r})$, but on all 
the orbitals explicitly, the HF equation is computationally more expensive in comparison with the DFT analogues.
The Fock matrix $\bold{F} = \bold{F}(\Phi)$ is defined as
$$
F_{\alpha\beta} = \int_{\mathbb{R}^3} \phi_\alpha H(\Phi) \phi_\beta \, d\bold{r}, \quad \alpha,\beta = \overline{1, N}.
$$
It is diagonal when $\Phi$ is the exact solution of \eqref{hf:general_form}.
Total energy can be calculated as
$$
E = 2\sum_{i=1}^N \lambda_i - \frac{1}{2} \iint_{\mathbb{R}^6} \frac{ \rho(\bold{r}) \rho(\bold{r'})}{| \bold{r} - \bold{r'}|} \, d\bold{r} d\bold{r'} + 
 \iint_{\mathbb{R}^6} \frac{|\tau(\bold{r}, \bold{r'})|^2 }{| \bold{r} - \bold{r'}|} \, d\bold{r} d\bold{r'},
+ E_{nn}
$$
for the Hatree-Fock model and
$$
E = 2\sum_{i=1}^N \lambda_i - \frac{1}{2} \iint_{\mathbb{R}^6} \frac{ \rho(\bold{r}) \rho(\bold{r'})}{| \bold{r} - \bold{r'}|} \, d\bold{r} d\bold{r'} + E_{xc}(\rho) + E_{nn},
$$
for the Kohn-Sham model, where
$$
E_{nn} = \sum_{i=1}^N \sum_{j>i}^N \frac{Z_i Z_j}{|\bold{R}_i - \bold{R}_j|},
$$
describes the repulsion between the nuclei.

\section{Tensor decomposition approach}\label{sec:tensors}

Orbitals $\phi_{\alpha}, \alpha=1, \ldots, N$ depend on three variables. Discretization on a rectangular $n \times n \times n$ mesh \footnote{the mode sizes can be in general 
different, but for simplicity we always assume that they are equal to $n$.} gives a three-dimensional tensor of coefficients. 
The storage complexity grows cubically in $n$. 
To get linear complexity $\mathcal{O}(n)$ we will use Tucker approximation \cite{tucker-factor-1963,lathauwer-svd-2000,khor-tuckertype-2007} of the all the tensors arising in the computations.

Tensor $A=\{a_{ijk}\}_{i,j,k=1}^n$ is said to be in the Tucker format if it is represented as
\begin{equation}\label{tucker}
a_{ijk} = \sum_{\alpha=1}^{r_1} \sum_{\beta=1}^{r_2} \sum_{\gamma=1}^{r_3} 
g_{\alpha\beta\gamma} u_{i\alpha} v_{j\beta} w_{k\gamma},
\end{equation}
where the minimal possible numbers $r_1, r_2, r_3$ required to
represent $A$ in the form \eqref{tucker} are called Tucker ranks,
matrices
$U = \{ u_{i\alpha}\}_{i=1,\alpha =1}^{n, r_1}$, 
$V = \{ v_{j\beta}\}_{j=1,\beta =1}^{n, r_2}$, 
$W = \{ w_{k\gamma}\}_{k=1,\gamma =1}^{n, r_3}$
are called factors, 
$G = \{g_{ijk}\}_{i,j,k=1}^n$
is called core of the decomposition.

Decomposition \eqref{tucker} contains $r^3 + 3nr$ parameters, compared with $n^3$ elements of the whole tensor.
In practice, instead of the decomposition we use the approximation with some accuracy~$\epsilon$.

Thus, each orbital on the grid is approximated in the Tucker format
$$
\phi_p (x_{i} ,y_{j}, z_{k}) \approx \sum_{\alpha=1}^{r_1(p)} \sum_{\beta=1}^{r_2(p)} \sum_{\gamma=1}^{r_3(p)} 
g^p_{\alpha\beta\gamma} u^p_{i\alpha} v^p_{j\beta} w^p_{k\gamma} + \mathcal{O}(\epsilon) ,\quad p = \overline{1, N} \quad i,j,k = \overline{1, n}
$$
Our main assumption is that the Tucker ranks are small and grow only logarithmically both in $n$ and $\varepsilon^{-1}$.
This has been verified before in \cite{mpi-chem3d-2009}. 
Once the approximation is computed, we can check the accuracy of the computations by computing relevant physical quantities and comparing it with the results, obtained by other methods.
\section{Iterative method}\label{sec:iterative}
The standard iterative technique to solve HF/KS equations is the self-consistent field iteration (SCF) \cite{jens-introqchem-2013}:
$$
H^{(k)}\ \phi^{(k+1)}_i = \lambda^{(k+1)}_i \phi^{(k+1)}_i, \quad i=\overline{1, N}. 
$$
This is not easy to implement in the Tucker format, which is a non-linear parametrization of the 3D-tensor, thus the original convex minimization problem 
is replaced by a non-convex one. 
Hence, SCF iteration is not much simpler that the original problem.
Therefore, we use a more convenient for our purposes Block Green iteration for Lippman-Schwinger form of \eqref{hf:general_form} as it can be easily done within tensor arithmetics.
In next subsection we describe the Block Green iteration technique with rotation of orbitals and present formulas to calculate the Fock matrix without direct Laplacian operation.
\subsection{The Block Green iteration}
Let us rewrite \eqref{hf:general_form} as follows
$$
\phi_i = -2 (-\Delta - 2 \lambda_i)^{-1} V\phi_i.
$$
The action of $ (-\Delta - 2 \lambda_i)^{-1}$ can be written as a convolution with the Yukawa potential
\begin{equation}\label{yukawa}
(-\Delta - 2 \lambda_i)^{-1} V\phi (\bold{r}) \equiv { \int_{\mathbb{R}^3}}\frac{e^{- \sqrt{-2\lambda_i}\ \|\bold{r} - \bold{r'}\|}}{4\pi \|\bold{r} - \bold{r'}\|}\ V\phi_i (\bold{r'})\ d\bold{r'}.
\end{equation}
We found that direct application of \eqref{yukawa} leads to numerical difficulties, which will be discussed later.
In this case much more efficient is to solve the screened Poisson equation directly using finite difference method.

For simplicity we start our description of the iteration for the system with one orbital $\phi_1 \equiv \phi$.
In this case the $k$-th step of the Green iteration is
$$
\tilde\phi = 2\left(-\Delta - 2 \lambda^{(k)}\right)^{-1} V^{(k)}\phi^{(k)}, \quad \phi^{(k+1)} = {\tilde\phi}/{\|\tilde\phi\|}.
$$
The energy is calculated as 
\begin{equation}\label{update1}
\lambda^{(k+1)} = \left( H^{(k+1)}  \phi^{(k+1)},\ \phi^{(k+1)} \right).
\end{equation}
The Hamiltonian $H^{(k)} \equiv H\left(\phi^{(k)}\right)$ contains the Laplacian $\Delta$. The approximation error in the Tucker format is controlled only 
in the discrete $L_2$ norm, thus discrete differentiation will amplify the error.
To avoid it \eqref{update1} can be rewritten as
\begin{equation}\label{update2}
\lambda^{(k+1)}  = \lambda^{(k)} + \frac{\left( V^{(k+1)}\tilde\phi - V^{(k)}\phi^{(k)},\ \tilde \phi\right)}{\left(\tilde\phi,\ \tilde\phi\right)}.
\end{equation}
The generalization to the case of molecules with more than one orbital is as follows.
We start by modifying each orbital separatedly
\begin{equation}\label{helm-iter}
\hat\phi_i =2\, (-\Delta - 2 \lambda_i^{(k)})^{-1} V^{(k)}\phi_i^{(k)},
\end{equation}
and orthogonalize $\widehat \Phi = \left(\hat\phi_1, \dots, \hat\phi_{N} \right)$ by calculating the Cholesky decomposition of the Gram matrix 
\begin{equation}\label{cholesky}
\int_{\mathbb{R}^3} \widehat \Phi^T \widehat \Phi \, d\bold{r}= LL^T, \quad \widetilde \Phi = \widehat \Phi L^{-T}.
\end{equation}
Then we calculate the $N\times N$ Fock matrix 
\begin{equation}\label{fock}
F =  \int_{\mathbb{R}^3} \widetilde\Phi^T  H^{(k+1)}\ \widetilde\Phi \, d\bold{r},
\end{equation}
diagonalize $F$
\begin{equation}\label{alg1}
 \Lambda^{(k+1)} = S^{-1} F S, 
\end{equation}
and finally rotate the orbitals by $S$
\begin{equation}\label{algfin}
\Phi^{(k+1)} = \widetilde \Phi S.
\end{equation}
The new values of orbital energies are the diagonal parts of $\Lambda^{(k+1)}$.
Steps of the iterative process are summarized in Algorithm \ref{alg-iter}
\begin{algorithm}[H] 
\begin{algorithmic}[1]
\caption{Block-Green iteration}\label{alg-iter}
\STATE Calculate $\widehat \Phi = (\hat \phi_1, \dots, \hat \phi_N)$: $\hat\phi_i =2\, (-\Delta - 2 \lambda_i^{(k)})^{-1} V^{(k)}\phi_i^{(k)}$.
\STATE Orthogonalize $\widehat \Phi$:  $\widetilde \Phi = \widehat \Phi L^{-T}$, where $L$: $LL^T =  \int_{\mathbb{R}^3} \widehat \Phi^T \widehat \Phi \, d\bold{r}$.
\STATE Calculate the Fock matrix $F =  \int_{\mathbb{R}^3} \widetilde\Phi^T  H^{(k+1)}\ \widetilde\Phi \, d\bold{r}$ via derivative-free formula in Statement \ref{statement-fock-update}.
\STATE Find new orbital energies by diagonalizing $F$: $ \Lambda^{(k+1)} = S^{-1} F S$
\STATE Find new orbitals: $\Phi^{(k+1)} = \widetilde \Phi S$
\end{algorithmic}
\end{algorithm}
\subsection{Derivative-free Fock matrix computation}
To compute the Fock matrix \eqref{fock} one need to find Laplacian of orbital functions. 
Since all computations in tensor formats are done approximately, differentiation will lead to the accuracy loss.
Therefore, we present the following derivative-free formula for the Fock matrix calculation

\newtheorem{Th}{Statement}
\begin{Th}\label{statement-fock-update}
In Algorithm  \ref{alg-iter} the Fock matrix \eqref{fock}  can be written in the form without derivatives
\begin{equation}\label{fock-update}
F = \int_{\mathbb{R}^3} \left( \widetilde\Phi^T  V^{(k+1)} \widetilde\Phi - \widetilde\Phi^T  V^{(k)} \Phi^{(k)}L^{-T} + L^{-1}\,\widehat\Phi^T\widehat\Phi\Lambda^{(k)} L^{-T} \right) \, d\bold{r},
\end{equation}
where $V^{(k)} \equiv V\left( \Phi^{(k)}\right)$.
\end{Th}

\begin{proof}
From \eqref{fock} we have
\begin{equation}\label{fock:proof}
F = L^{-1} \left(\int_{\mathbb{R}^3} \widehat\Phi^T  H^{(k+1)}\ \widehat \Phi\, d\bold{r} \right)  L^{-T},
\end{equation}
Matrix $$\int_{\mathbb{R}^3} \widehat\Phi^T  H^{(k+1)}\ \widehat \Phi\, d\bold{r}$$ can be written elementwise in the following way
\begin{align}
\int_{\mathbb{R}^3} \hat\phi_\alpha  H^{(k+1)} \hat\phi_\beta\, dx =
\int_{\mathbb{R}^3} \hat\phi_\alpha \left( -\frac{1}{2}\Delta + V^{(k+1)}  \right) \hat\phi_\beta   \, dx =\notag \\
\int_{\mathbb{R}^3} \hat\phi_\alpha  V^{(k+1)}  \hat\phi_\beta   \, dx + 
\int_{\mathbb{R}^3} \hat\phi_\alpha \left( -\frac{1}{2}\Delta \pm \lambda^{(k)}_\beta \right) \hat\phi_\beta   \, dx = \notag \\
\int_{\mathbb{R}^3} \hat\phi_\alpha  \left( V^{(k)}  + \lambda^{(k)}_\beta  \right) \hat\phi_\beta \, dx + 
\int_{\mathbb{R}^3} \hat\phi_\alpha \left( -\frac{1}{2}\Delta - \lambda^{(k)}_\beta \right) \hat\phi_\beta   \, dx \notag
\end{align}
Taking into account that  $$\hat\phi_\beta = \left(-\frac{1}{2}\Delta - \lambda_\beta^{(k)}\right)^{-1} V^{(k)}\phi_\beta^{(k)} $$ we have
\begin{align}\notag
\int_{\mathbb{R}^3} \hat\phi_\alpha \left( -\frac{1}{2}\Delta - \lambda^{(k)}_\beta \right) \widetilde\phi_\beta   \, dx = 
\int_{\mathbb{R}^3} \hat\phi_\alpha \left( -\frac{1}{2}\Delta - \lambda^{(k)}_\beta \right)  \left(-\frac{1}{2}\Delta - \lambda_\beta^{(k)}\right)^{-1} V^{(k)}\phi_\beta^{(k)}  \, dx = 
\int_{\mathbb{R}^3} \hat\phi_\alpha  V^{(k)}\phi_\beta^{(k)}  \, dx. \quad
\end{align}
Finally, substituting the obtained expressions into \eqref{fock:proof} we get \eqref{fock-update}.
\end{proof}
\subsection{Density mixing}

In order to accelerate convergence we used the density mixing scheme, which is typical used for DFT calculations \cite{jens-introqchem-2013}.
Here we provide a brief description of the scheme for the sake of completeness.
Our scheme is a fix point iteration which can be formally written as $\rho^{(k)} = G(\rho^{(k-1)})$. In the density mixing scheme, 
the next density $\rho^{(k+1)} $ is represented as  a linear combination of the results of previous iterations
$$
\rho^{(k+1)} = (1-\beta_k) \sum_{j = 1}^m \alpha_j  \rho^{(k - m + j)} + \beta_k\sum_{j = 1}^m \alpha_j G\left( \rho^{(k - m + j)} \right),
$$
where the coefficients $\alpha = (\alpha_1, \dots, \alpha_m)$ are the solution of the minimization problem
$$
\alpha = \arg \min_{\tilde\alpha_1, \dots, \tilde\alpha_m} \| \sum_{j = 0}^{m} \tilde\alpha_j \left[ \rho^{(k - m + j)} - G\left( \rho^{(k - m + j)} \right) \right] \| 
$$
under the constraints
$$
\sum_{j = 0}^{m} \tilde\alpha_j  = 1.
$$
\section{Discretization}\label{sec:discr}
It is known that the orbitals decay exponentially at infinity \cite{ishida-decay-1992}:
\begin{equation} \label{decay}
\phi_i(\bold{r}) = \mathcal{O} \left(  e^{-\sqrt{-2\lambda_i} \, \|\bold{r}\|}\right) , \quad \|\bold{r}\|\to+\infty.
\end{equation}
As a result, multiplications and linear combinations of the orbitals also decay exponentially as $\|\bold{r}\| \to +\infty$. 
Hence, we replace the whole $\mathbb{R}^3$ with a finite box $\Omega = [-L, L]^3$, where $L$ depends on the chosen accuracy and estimation of the smallest orbital energy $\lambda_{N}$. 
We will numerically investigate how the choice of $L$ in the experiments section.

In $\Omega$ we introduce a uniform tensor-product grid $\omega^h = \omega_1^{h_1} \times \omega_2^{h_2} \times\omega_3^{h_3}$ with $h = 2L_i/n_i$, where $\omega_i^h = \{ -L_i + kh_i: \, k=0, ..., n \}$, $i=1,2,3$. We use uniform grids to obtain structured operators and significantly decrease computation complexity. Low-order approximation accuracy over uniform grids will be improved later by using extrapolation. 
For simplicity we use grids with $h_1 = h_2 = h_3$ and $L_1 = L_2 = L_3$.
Algorithm \eqref{helm-iter}--\eqref{algfin} requires the computation of $V\Phi$ at each step, and the 3D-convolution have to be computed 
\begin{equation} \label{conv-newt}
w(\bold{r}) \equiv
\int_{\mathbb{R}^3} \frac{f(\bold{r'})}{\|\bold{r} - \bold{r'}\|} \, d\bold{r'}, \quad \bold{r}\in \mathbb{R}^3. 
\end{equation}
\begin{remark}
One can use the fact that problem of finding $w(\bold{r})$ is equivalent to find the solution of $-\Delta w =4\pi f$ and it seems better to solve less expensive Poisson equation instead of direct computation of the convolution.
Nevertheless, $w(\bold{r}) = \mathcal{O}(1/\|\bold{r}\|)$ when $\|\bold{r}\|\to\infty$. 
Boundary conditions are unknown and, hence, to get the accuracy $\epsilon$ one have to choose domain size $L = \mathcal{O}\left( 1/\epsilon\right)$. 
\end{remark}
On the other hand we can compute convolution with the Yukawa kernel (operator $(-\Delta -2\lambda)^{-1}$ in \eqref{helm-iter}) by directly inverting the shifted Laplacian.
This can be done efficiently due to the exponential decay of orbitals~\eqref{decay}, so it is safe to enforce Dirichlet boundary conditions.
Moreover, the direct convolution with the Yukawa kernel discretized on a uniform grid leads to the non-symmetry of the Fock matrix and this breaks down the convergence.
The point here is that the numerical analogue of the derivative-free formula \eqref{fock-update} does not hold unless the  action of $(-\Delta -2\lambda)^{-1}$ is computed by solving screened Poisson equation.

Thus, to discretize $(-\Delta -2\lambda)^{-1}$ we use standard 7-points finite difference scheme.
To discretize convolutions \eqref{conv-newt} we use Galerkin scheme and piecewise-constant basis functions $\phi_{\bold{i}}$ with support on $\Omega_{\bold{i}}$, where $\bold{i} \in \mathcal{I} \equiv \{0,\dots, n-1\}^3$ and $\Omega_{\bold{i}}$ are $h^3$ cubes centered in $r_{\bold{i}} = (x_i, y_j, z_k)\in \omega^h$. 
So,
\begin{equation}\label{conv-discretized}
w(\bold{r}_i)\equiv w_\bold{i} \approx \sum_{\bold{j}\in\mathcal{I}} f_\bold{j}\,  q_{\bold{i}-\bold{j}},
\end{equation}
where
\begin{equation} \label{galerk}
	q_{\bold{i}-\bold{j}} =  \int_{\mathbb{R}^3}\int_{\mathbb{R}^3} \frac{\phi_{\bold{i}}(\bold{r}) \phi_{\bold{j}}(\bold{r'})}{\|\bold{r} - \bold{r'}\|} \, d\bold{r} d\bold{r'}, \quad f_{\bold{i}} = \int_{\mathbb{R}^3} f(\bold{r}) \phi_{\bold{i}}(\bold{r})  \, d\bold{r}.
\end{equation}
The total approximation error of the discretization scheme is $\mathcal{O}(h^2)$. 
To get high-order discretization schemes translation-invariant basis functions of higher order can be used $\phi_\bold{i}(\bold{r}) = \phi(\bold{r} - \bold{r}_\bold{i})$, 
where $\phi(\bold{r})$ is a suitable piecewise-polynomial function.

\section{Operations in the Tucker format}\label{sec:operations}
In order to implement the iterative algorithm in the Tucker format, we need several basic operations. 
Linear operations (addition, multiplication by number) can be done straightforwardly \cite{ost-latensor-2009}.  
Calculation of scalar products and norms can also be efficiently implemented within the tensor formats. After such operations the ranks typically increase, 
but can be reduced thanks to the SVD-based rounding procedure. 
\subsection{Elementwise operations and cross approximation}
Despite linear complexity in the mode size, some operations such as elementwise products and convolutions may have strong rank dependence: $\mathcal{O}(r^8)$ floating point operations.
To decrease this complexity we use the so-called \emph{cross approximation method} \cite{tee-cross-2000, bebe-2000, ost-tucker-2008}.
This method finds the decomposition using only few of elements of the approximated array.
Every next iteration of the cross approximation procedure chooses adaptively new elements to be computed untill the stopping criterea holds.
In the implementation we use Schur-Cross3D algorithm proposed in \cite{ro-crossconv-2015} that has $\mathcal{O}(r^3 + nr)$ complexity. 

The cross approximation algorithm is also the key technique to approximate nonlinear elementwise functions of density arising in the calculation of  exchange-correlation potential $V_{xc}$.
We note, that $V_{xc}(\rho)$ itself  does not have low-rank structure, however, $V_{xc}(\rho)\,\phi_i$ is of low rank. 
\subsection{Convolution}
Computation of the convolution is one of the most expensive parts of \eqref{hf:general_form} with strong dependence on ranks of input tensors.
We calculate it by the \emph{Cross-Conv} algorithm proposed in \cite{ro-crossconv-2015}. This algorithm has  $\mathcal{O}(r^4 + nr^2)$ complexity, which is the fastest known for practically interesting mode sizes $n$ up to $n\sim 2^{15}$.
However, for fine meshes it might be a good idea to use quantize tensor train (QTT) convolution algorithm \cite{khkaz-conv-2013} can be used. This algorithm 
computes the convolution of two tensors given in the Tucker format. The orbitals are already there, and it is well-known that the kernel function ${1}/{r}$ can 
be approximated in the Tucker format with $r = \mathcal{O}(\log n \log^{-1} \epsilon)$.
This trick with decomposition into sum of exponents was used to calculate external potential \eqref{ext-pot}.
\subsection{Fast Laplacian inversion}
The Laplacian equation is solved using the algorithm proposed in \cite{om-stockes-2010}. 
Given the right-hand side in the Tucker format, the solution after the sine transform (which is separable) reduces to the computation of 
$$
   \frac{a_{ijk}}{\eta_i + \eta_j + \eta_k - \mu} ,
$$
where $\eta_i, i = 1, \ldots, n$ are the eigenvalues of the one-dimensional Laplace operator. 
Such computation can be straightforwardly done using the cross-approximation technique.

\section{Complexity estimates}\label{sec:complexity}

For simplicity we provide an upper bound of complexity for $r$ be the largest rank of all orbitals.
We also assume that the calculations are done on $n\times n \times n$ grid.
Precomputations before the iterative process include  calculation of the external potential \eqref{ext-pot} and the Galerkin tensor of the convolution kernel \eqref{galerk}.
The computation of the external potential requires $N_{nuc}\cdot \mathcal{O} (r^3 + 3nr)$.
Precomputation of the Galerkin tensor requires $\mathcal{O} (r^3 + 3nr)$ operations.
It does not depend on number of nuclei but has a bigger constant due to one-dimensional integration.

Let us now estimate the complexity of one iteration.
We denote by $K_{cross} = \mathcal{O}(r^4 + nr^2)$ the complexity of the cross approximation calculated by Schur-Cross3D and by $K_{conv} = \mathcal{O}(r^4 + nr^2 + rn \log n)$ the complexity of one convolution operation \cite{ro-crossconv-2015}. 
We note that the constant hidden in $\mathcal{O}(\cdot)$ is of the order of unity.
Calculation of $V^{(k)} \phi_i^{(k)}$ in \eqref{helm-iter} is done as follows.
First of all we calculate $V_{coul}$ which is convolution operation that has $K_{conv}$ complexity.
The Coulomb potential is calculated only once on each iteration.
In case of KS equations we run cross approximation procedure to calculate $V^{(k)} \phi_i^{(k)}$ since the $V_{xc}$ does not have low-rank structure.
For the HF equations we need additionally to calculate the exchange potential, which requires $N^2$ convolution operations and $N^2$ elementwise products.
Thus, the step has $N K_{cross}$ complexity for KS and $N^2 K_{cross}$ for the HF equations.
The application of $(-\Delta - \mu)^{-1}$ is also of $N K_{conv}$ complexity but with approximately twice smaller constant in $K_{conv}$ as one of the input tensors has fixed rank 3.

Step \eqref{cholesky} requires $N^2$ scalar product evaluations, so it has $N^2  K_{cross}$ complexity as elementwise products are done by the cross approximation. Complexity $N^2  K_{cross}$ is much larger than the Cholesky factorization as we assume that $N \ll r^4$.
The Fock matrix computation \eqref{fock-update} consists of $V^{(k)}\widetilde \Phi$ which is of same complexity as $V^{(k)} \Phi^{(k)}$ and of scalar products.
Thus,  we get the overall complexity of one iteration to be $N^2 \cdot \mathcal{O}(r^4 + nr^2 + rn\log n)$.

\section{Numerical experiments} \label{sec:num}
The prototype is implemented in Python. 
The implementation of HF and KS solvers can be found at \url{https://github.com/rakhuba/tensorchem}. 
The toolbox with arithmetics in the Tucker format can be found at \url{https://github.com/rakhuba/tucker3d}. 
For the basic linear algebra tasks the MKL library is used. Python and MKL are from the Anaconda
Python Distribution \url{https://store.continuum.io/cshop/anaconda/}. 
Python version is 2.7.9. MKL version is 11.1-1. Tests were performed on 4 Intel Core i7 2.6
GHz processor with 8GB of RAM. However, only 2 threads were used (this is default
number of threads for MKL). We would like to emphasize that implementation of the
whole algorithm is in Python and time performance can be improved by
implementing the most time-consuming parts of it in C or Fortran languages.


\subsection{The box size}\label{sec-num-box}
\begin{figure}[h!]\label{pic-box}
\begin{center} 
\resizebox{0.7\textwidth}{!}{\input{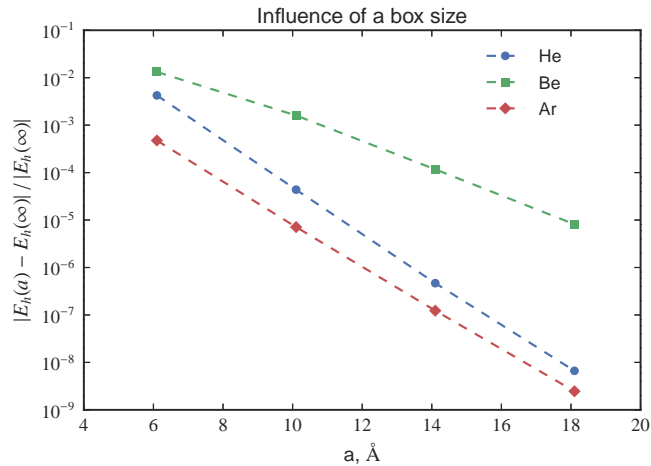}}
\end{center}
\caption{Dependence of relative error of $E_h(a)$ on a box size $a$ for fixed $h=0.1$ \AA,\ $\epsilon = 10^{-9}$. To calculate $E_h(\infty)$ the box of size $a = 50$ \AA\  was used.} 
\label{pic:3d}
 \end{figure}

Here we illustrate how the results depend on the simulation box size.
The orbitals decay exponentatially as $\exp (-\sqrt{-2\lambda_{\text{HOMO}}}\|x\|)$, $\|x\|\to\infty$. 
Therefore, one can expect that the error introduced by the finite size of the box has exponential decay with the box size.
On Figure \ref{pic-box} illustrates this guess.
On this figure the dependence of relative error of $E_h(a)$ with respect to the $E_h(\infty)$ as a function of a box size $a$ is presented.
All calculations were done for the fixed grid step $h=0.1$ \AA\ and accuracy $\epsilon = 10^{-9}$.
$E_h(\infty)$ was estimated at $a = 50$~\AA.

\subsection{Extrapolation}
As the proposed approach is fully grid-based, a good idea is to start the method on a coarse grid with, say $N = 128$ and then use it as an initial guess on finer grids (grids with $2^{k}$ grid points are used). 
Due to the linear comlexity of the proposed algorithm the total time for extrapolation is not larger than twice the time required for the computation on the finest grid:
\begin{equation}\label{exrapolation-time}
t_\text{extrapolated} = t_N + t_{N/2} + \dots + t_{N/2^{l}} = C N + CN/2 + \dots + CN/2^{l}  < 2CN = 2t_\text{N},
\end{equation}
where constant $C$ does not depend on $N$.
This makes the extrapolation a very useful and computationally efficient part of the whole algorithm. 

\begin{figure}[h!]\label{pic-extrapolation}
\begin{center} 
\resizebox{0.7\textwidth}{!}{\input{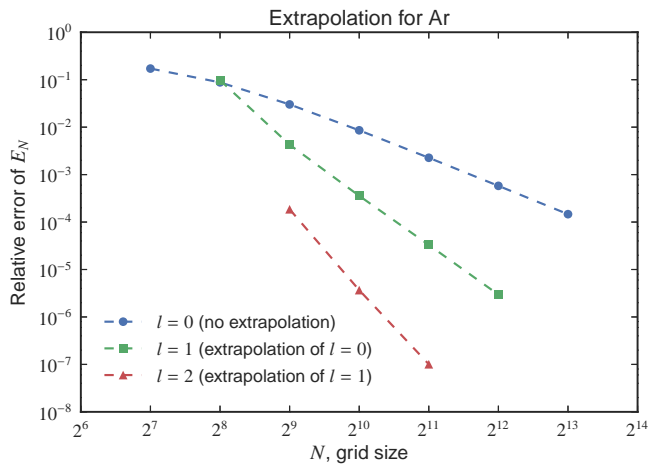}}
\end{center}
\caption{Extrapolation of the full energy as a function of grid size for Argon. Error is calculated with respect to highly accurate results from \cite{thakkar-atoms-1995}.} 
\label{pic:3d}
 \end{figure}

Numerical results showed that on coarse grids the order of the approximation is not exactly of the second order. 
Hence, we used the so-called Aitken's delta squared process \cite{joyce-extrap-1971}, which is equivalent to the Richardson extrapolation for the exact second order. The Aitken's process accelerates the convergence of the sequence $E_n$ such that the new sequence $E'_n$
\begin{equation}\label{aitken}
 E'_n = E_{n+2} - \frac{(E_{n+2} - E_{n+1})^2}{E_{n+2} - 2E_{n+1} + E_{n}},
\end{equation}
converges faster than $E_n$ when $n$ goes to infinity. The same trick can be applied to $E'_n$.
The results for extrapolation are shown on Figure~\ref{pic-extrapolation}.

\subsection{Hartree-Fock calculations on atoms}
First we present the HF calculations of closed-shell atoms He, Be, Ne, Ar.
Table \ref{tab-atoms} represents total and  highest occupied molecular orbital (HOMO)  energies.
These results are compared with highly accurate results for atoms \cite{thakkar-atoms-1995}.
Our grid-based results are obtained with the relative accuracy $\epsilon=10^{-7}$. 
We use the relative error instead of absolute in order to obtain accurate HOMO values and at the same time to be faster in the calculation of lower orbitals.
\begin{table}\label{tab-he}
\begin{center}
\begin{tabular}{c|cccccc}
\hline
$\epsilon$ & $10^{-3}$ & $10^{-5}$ & $10^{-7}$ & $10^{-9}$ & $10^{-11}$ &   \\
Total energy & -{2.861}5& -{2.86167}8  & -{2.861680}1  & -{2.86168000}0 &-{2.8616799959}3 \\
\hline
\end{tabular}
\end{center}
\caption{Helium total energy for different values of relative error $\epsilon$. The energy was extrapolated on a sequence of grids: from $128^3$ to $8196^3$. Bold denotes correct numbers. The total energy in  \cite{thakkar-atoms-1995} is -2.861 679 996.}
\end{table}
However, absolute accuracy can also be used.
Extrapolation is done on a range of grid sizes: from $128^3$ to $8192^3$. 

Results in Table \ref{tab-atoms} show systematic convergence of the total and HOMO energies. 
Note, that for all atoms under consideration extrapolation with grid sizes up to $N^3 = 8192^3$ grid points is  enough to get accuracy $\epsilon=10^{-7}$.
Without extrapolation, much larger grids are needed.
\begin{table}\label{tab-atoms}
\begin{center}
\begin{tabular}{c|ccccc}
\hline
Atom & Method & Total energy & HOMO energy  \\
\hline 
He & $N=8192$ & -2.861 670  & -0.917 950  \\
& Extrapolated & -2.861 680  & -0.917 955  \\
& \cite{thakkar-atoms-1995} & -2.861 680  & -0.917 956 \\
& aug-cc-pVQZ & -2.861 539  &  -0.917 915 \\
& aug-cc-pV5Z &  -2.861 635  & -0.917 935 \\
\hline
Be & $N=8192$ & -14.572 256  & -0.309 263  \\
& Extrapolated & -14.573 023  & -0.309 269 \\
& \cite{thakkar-atoms-1995} & -14.573 023  & -0.309 270 \\
& aug-cc-pVQZ & -14.572 976   &  -0.309 269 \\
& aug-cc-pV5Z &  -14.573 011  & -0.309 270\\
\hline
Ne & $N=8192$ & -128.518 74 & -0.850 523 \\
& Extrapolated & -128.547 08 & -0.850 410  \\
& \cite{thakkar-atoms-1995} & -128.547 09 & -0.850 410 \\
& aug-cc-pVQZ & -128.544 69 & -0.850 210 \\
& aug-cc-pV5Z &  -128.546 87 & -0.850 391\\
\hline
Ar & $N=8192$ & -526.740 5 & -0.591 024 \\
& Extrapolated & -526.817 4 & -0.590 017  \\
& \cite{thakkar-atoms-1995} & -526.817 5 & -0.591 017 \\
& aug-cc-pVQZ & -526.816 9 & -0.591 013 \\
& aug-cc-pV5Z &  -526.817 3 &  -0.591 011 \\
\end{tabular}
\end{center}
\caption{Total and HOMO energies for different atoms, $\epsilon = 10^{-7}$. Extrapolated calculations were done on a sequence of grids: from $128^3$ to $8196^3$. aug-cc-pVQZ and  aug-cc-pV5Z calculations were done by NWCHEM  \cite{nwchem-2010}.}
\end{table}

\subsection{LDA calculations on molecules}
We perform the LDA calculations for several molecules. Molecular geometries were taken from the NIST database \cite{nist-2001}.
In Table \ref{tab-molecules} the resulting total and HOMO energies are compared with the LDA computations with the aug-cc-pVXZ (X$=$Q, 5) basis sets. 
The latter we done using NWCHEM program package \cite{nwchem-2010}.
In all experiments the relative error for all orbital functions was set to be $\epsilon=10^{-7}$.
The box size was chosen adaptively to the HOMO energy.

\begin{figure}[h!]\label{pic-convergence}
\begin{center} 
\resizebox{0.7\textwidth}{!}{\input{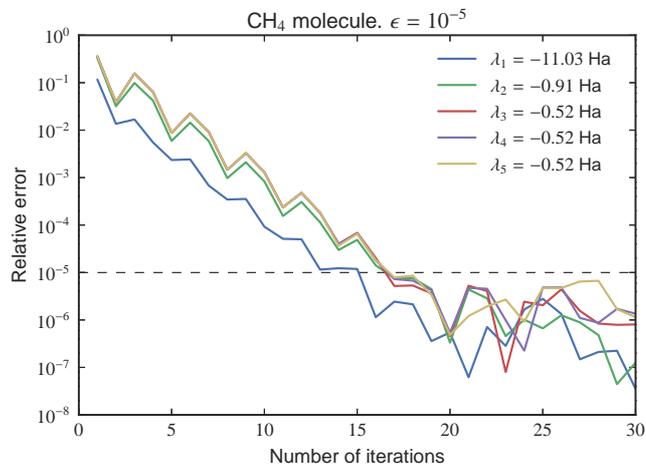}}
\end{center}
\caption{Relative errors of each orbital energy as a function of number of iterations ($\epsilon=10^{-5}$).} 
\label{pic:3d}
 \end{figure}

We note that the extrapolated energy is smaller than for the basis aug-cc-pV5Z approximately $10^{-4}$ hartree both for the total energy and for the HOMO energy, thus our method has less
basis set error.
\begin{table}\label{tab-molecules}
\begin{center}
\begin{tabular}{c|ccccc}
\hline
Molecule & Method & Total energy & HOMO energy  \\
\hline 
H$_2$ & $N=8192$ & -1.137 392 3 & -0.378 667  \\
& Extrapolated & -1.137 392 8 & -0.378 668 \\
& aug-cc-pVQZ & -1.137 249 9 & -0.378 649 \\
& aug-cc-pV5Z & -1.137 374 8&  -0.378 665\\
\hline 
CH$_4$ & $N=8192$ & -40.116 452  &  -0.348 989 \\
& Extrapolated & -40.119 829  & -0.348 984 \\
& aug-cc-pVQZ &  -40.118 644 & -0.348 964  \\
& aug-cc-pV5Z &  -40.119 299 & -0.348 982 \\
\hline 
C$_2$H$_6$ & $N=8192$ & -79.069 964 & -0.299 763  \\
& Extrapolated & -79.075 142 & -0.299 761\\
& aug-cc-pVQZ & -79.070 784 & -0.299 724\\
& aug-cc-pV5Z & -79.072 763& -0.299 762 \\
\end{tabular}
\end{center}
\caption{Total and HOMO energies for different molecules,  relative error $\epsilon = 10^{-7}$. Extrapolated calculations were done on a sequence of grids: from $128^3$ to $8196^3$. aug-cc-pVQZ and  aug-cc-pV5Z calculations were done by NWCHEM  \cite{nwchem-2010}.}
\end{table}
Figure \ref{pic-convergence} illustrates rapid convergence of all five occupied orbitals with respect to the number of iterations for $\epsilon = 10^{-5}$.
As anticipated the convergence stopped below $\epsilon$.
We note that in all experiments orbitals have converged to the $\epsilon = 10^{-7}$ within $30$ iterations.

\subsection{LDA calculations on cluster of atoms}
Here we present calculations on a test system of hydrogen atoms, which was firstly calculated in \cite{vekh-blackbox-2013} as an artificial example.
The system represents finite cluster of hydrogen atoms with unit cell having primitive cubic structure.
The distance between atoms is set to be $d = 2$~\AA.
As in section with LDA calculations on molecules in all experiments we observed systematical convergence of the total energy with $\epsilon$ precision.
Table \ref{tab-cluster-ranks} illustrates that ranks of orbitals grow sublinearly with the system size.
Thus, we expect the proposed algorithm to be  particularly efficient for systems with regular location of atoms.
The time of one iteration on $N=1024$ and $\epsilon = 10^{-5}$ is $9$ sec for $H_{3\times 2 \times 2}$ and $68$ sec for $H_{8\times 2 \times 2}$, which as anticipated scales approximately quadratically with the size of the system: $(8/3)^2 \approx 7.1$ and $68/9 \approx 7.5$.

\begin{table}\label{tab-cluster-ranks}
\begin{center}
\begin{tabular}{c|ccc}
\hline
Cluster & min ranks orbital & max ranks orbital\\
\hline 
H$_{1\times 2 \times 2}$ & $16 \times 16 \times 14$ & $17 \times 17 \times 15$ \\
H$_{3\times 2 \times 2}$ & $28 \times 28 \times 21$ & $35 \times 35 \times 20$ \\
H$_{8\times 2 \times 2}$ & $25 \times 25 \times 22$ & $36 \times 36 \times 26$ \\
\hline
H$_{1\times 2 \times 1}$ &  $13 \times 12 \times 12$ & $13 \times 12 \times 12$ \\
H$_{2\times 2 \times 1}$ &  $16 \times 16 \times 14$ & $17 \times 17 \times 15$ \\
H$_{5\times 2 \times 1}$ & $16 \times 18 \times 11$ & $16 \times 20 \times 13$ \\
H$_{9\times 2 \times 1}$ & $14 \times 18 \times 11$ & $17 \times 22 \times 13$ \\
H$_{16\times 2 \times 1}$ & $16 \times 19 \times 11$ & $21 \times 27 \times 14$ \\
\hline
H$_{1\times 4 \times 1}$ & $13 \times 12 \times 13$ & $13 \times 13 \times 13$ \\
H$_{3\times 4 \times 1}$ & $19 \times 20 \times 13$ & $23 \times 23 \times 15$ \\
H$_{8\times 4 \times 1}$ & $23 \times 20 \times 12$ & $32 \times 27 \times 15$ \\
\end{tabular}
\end{center}
\caption{Ranks of orbitals with maximum  and minimum ranks on different clusters of hydrogen atoms.}
\end{table}

\section{Conclusions and future work}
Tensor approach to the solution of 3D problems in electronic computations allows us to reach desirable accuracy at moderate computational cost. 
This can be used, for example, in verification of other methods and in construction of accurate basis sets for other methods  \cite{lrs-dmrgcc-2014}.
The main problem with the presented approach is still in its scaling for large $N$: for example, one has to compute the approximation of all the products $\phi_i \phi_j$ using the cross method. This can be solved by the resolution of identity methods \cite{vaf-ri-1993}, but they are typically used for global basis sets. 
What is interesting, is that the density itself shows a very good low-rank structure; thus, it is very interesting to apply tensor decomposition methods for the 
orbital-free DFT functionals \cite{gav-ofdft-2007, kg-ofdft-2012}. 

\section{Acknowledgements}
We thank DrSci. Boris Khoromskij and Dr. Venera Khoromskaia for useful comments on the draft of the manuscript.

\bibliographystyle{elsarticle-num}
\bibliography{bibtex/tensor,bibtex/our,bibtex/dmrg,quant.bib}
\end{document}